\documentclass[twocolumn]{autart}    

\usepackage{graphicx}       

\usepackage{amsmath}
\usepackage{color}
\usepackage{amssymb}
\usepackage{dsfont}
\usepackage{enumerate}

\def\ao#1{{\color{black}#1}}

\def\rmj#1{{#1}}

\def\aor#1{{#1}}

\let\theoremstyle\relax
\usepackage{epsfig} 
\usepackage{graphicx} 
\usepackage{tikz, tkz-berge} 
\usetikzlibrary {positioning} 

\begin{document}



\makeatletter
\DeclareRobustCommand*\cal{\@fontswitch\relax\mathcal}
\makeatother
\newcommand{\cerny}{{\v Cern{\'y} }{}} 
\newcommand{\re}{{\mathbb R}} 
\newcommand{\setmat}{{\cal{M}}} 
\newcommand{\cm}{{\cal{M}}} 
\newcommand{\consmat}{ { \mbox{E} } }
\newcommand\mb[1]{\mbox{#1}}
\newcommand{\G}{\mathcal{G}}
\newcommand{\D}{\mathcal{D}}
\newcommand{\R}{\mathcal{R}}
\newcommand{\V}{\mathcal{V}}
\newcommand{\E}{\mathcal{E}}
\newcommand{\Ed}{\mathcal{E}_{d}}
\newcommand{\graph}{ \G = \left( \V , \E \right) }
\newcommand{\id}{\mathds{I}}
\newcommand{\ev}{\mathbf{v}}
\newcommand{\1}{\mathbf{1}}
\newcommand{\digraph}{\mathcal{D} = (\V, \Ed)}
\newcommand{\vectun}{{e }} 
\newcommand{\qedsymbol}{\rule{1ex}{1ex}}

\newtheorem{ex}{Example} 
\newtheorem{opq}{Open Question} 

\newtheorem{deff}[thm]{Definition} 
\newtheorem{assumption}{Assumption} 

\begin{frontmatter}

\title{On Primitivity of Sets of Matrices}

\author[Uclouvain]{V. Blondel}, 
\author[Uclouvain]{R. M. Jungers}, 
\author[UIUC]{A. Olshevsky} 

\address[Uclouvain]{Universit\'e catholique de Louvain, Belgium}  
\address[UIUC]{University of Illinois at Urbana-Champaign, USA}             

 \begin{abstract}A nonnegative matrix $A$ is called primitive if $A^k$ is positive for some integer $k>0$.  A generalization of this concept to \aor{finite} sets of matrices is as follows: a set of matrices $\setmat = \{A_1, A_2, \ldots, A_m \}$ is primitive if $A_{i_1} A_{i_2} \ldots A_{i_k}$ is positive  
for some indices $i_1, i_2, ..., i_k$. The concept of primitive sets of matrices \aor{comes up in a number of problems within the study of discrete-time switched systems}. In this paper, we analyze the computational complexity of deciding if a given set of matrices is primitive and we derive bounds on the length of the shortest positive product.  

We show that while primitivity is algorithmically decidable, unless $P=NP$ it is not possible to decide primitivity of a matrix set in polynomial time. Moreover, we show that the length of the shortest positive sequence can be superpolynomial in the dimension of the matrices. On the other hand, defining ${\cal P}$ to be 
the set of matrices with no zero rows or columns, we give a simple combinatorial proof \rmj{of a previously-known characterization of primitivity for matrices in ${\cal P}$ which can be} tested in polynomial time.  This latter observation is related to the well-known 1964 conjecture of \cerny    on synchronizing automata; in fact, any bound on the minimal length of a 
synchronizing word for \rmj{synchronizing} automata immediately translates into a bound on the length of the shortest positive product of a primitive set of matrices in ${\cal P}$. In particular, any primitive set of $n \times n$ matrices in ${\cal P}$ has a positive product of length $O(n^3)$. \end{abstract} 

\end{frontmatter}

\section{Introduction}

A $n \times n$  matrix $A$ \aor{which is entrywise nonnegative} is said to be primitive if every entry of $A^k $ is positive for some positive integer $k$.  It is well-known (see \cite{HJ_Matrix_Analysis_Book}, Corollary 8.5.9) that this is the case if and only if $A^{n^2-2n+2}>0$ so the primitivity of a matrix is easy to verify algorithmically.  A straightforward generalization of \aor{primitivity to finite} sets of matrices is the following \cite{PV}: a finite set of  $m$ nonnegative matrices $\setmat =\{A_1, A_2, \ldots, A_m \}$ is primitive if  $A_{i_1} A_{i_2} \ldots A_{i_k} $ \aor{is (entrywise) positive} for some indices $i_1, i_2, ..., i_k \in \{1, ..., m\}$.  

The property of primitivity of a set of matrices is important in several applications.  In particular, its presence enables one to use efficient algorithms for the computation of the {Lyapunov exponent} of a stochastic switching system (we refer the reader to \cite{liberzon-switching,shorten-siamreview07,jungers_lncis} for a general introduction to switching systems).  Given a finite set of matrices $\setmat \subset \re^{n \times n},$ one can define a stochastic switched system as:  \rmj{
\begin{equation} 
x_{k+1}=A_{i_k}x_k, \quad A_{i_k}\in \setmat,\label{eq:switched.system} 
\end{equation}} where for simplicity let us make the assumption that each $A_{i_k}$ is chosen randomly from the uniform distribution on $\setmat$. Such models are commonly used throughout stochastic control; for example, they are a common choice for modeling manufacturing systems with random component failures (see \cite[Chapter 1]{stochbook}). The \emph{Lyapunov exponent} of this system is defined by the following limit \rmj{(where $\E$ denotes the expectation): }
\begin{equation}\label{LE} 
\lambda \quad = \quad \lim_{k \to \infty} \ \frac1k \ \E \log \ \bigl\| \, A_{d_k}\cdots A_{d_1} \bigr\|\, . 
\end{equation} 
 The Lyapunov exponent characterizes the rate of growth of the switching system \rmj{with probability one}.  
%
While it is hard to compute in general \cite{tsitsiklis97lyapunov}, it turns out that {\em in the particular case of primitive sets of matrices, efficient algorithms are available.} \ao{We refer the reader to \cite{jungers-protasov-lyapunov_ecc,jungers-protasov-lyapunov,protasov-invariant-english,protasov-primitive} for the algorithms.}   

\ao{Secondly, the concept of primitivity is also related to the so-called consensus problem. Here the matrices in $\setmat$ are further taken to be stochastic matrices and the question is whether the recursion of Eq. (\ref{eq:switched.system}) \aor{almost surely} converges to $\alpha {\bf 1}$, i.e., to a multiple of the all-ones vector. In this case, we say that the iteration achieves consensus on {the value} $\alpha$. Such ``consensus iterations'' appear in a number of applications, and there  is now a considerable literature on the consensus problems providing necessary or sufficient conditions under various assumptions on the switching - we refer the reader to the classical and modern papers \cite{cons1}, \cite{cons2}, \cite{cons3}, \cite{cons4}, \cite{cons5}, \cite{cons6}, \cite{bl-ol} for examples of such conditions and discussions of applications.   }

\ao{\aor{The consensus problem naturally leads to the concept of primitivity when, as before, the matrices $A_{i_k}$ are chosen randomly}. Indeed, if the matrices $A_{i_k}$ are chosen from
the uniform distribution on $\setmat$ (or from any distribution whose support is $\setmat$)  and if we further stipulate that $\alpha$ should lie in the interior of the convex hull of the entries of \rmj{$x_0$} (which means every node has an influence on the final consensus value), then it is immediate that consensus on such an $\alpha$ is achieved if and only if $\setmat$ is primitive.} 

\aor{Finally, the problem of matrix primitivity is perhaps the simplest possible reachability problem for switched systems: given an unknown initial state  in the nonnegative orthant, can we choose at each step a matrix $A_{i_k}$ from the set of nonnegative matrices $\{A_1, \ldots, A_m\}$ so that the final state of Eq. (\ref{eq:switched.system}) is in the interior of the nonnegative orthant? As we show in this paper, even this simple and stylized reachability problem faces significant computational obstructions. }

In this paper, we study the problem of recognizing primitivity and related problems. Given a set of $n \times n$ nonnegative matrices $\setmat = \{A_1, \ldots, A_m\}$ one would like to determine, efficiently if possible, whether or not $\setmat$ is primitive. This is closely related to the problem of bounding the length of the shortest positive product of matrices from $\setmat$, which we denote by $l(\setmat)$. Indeed, an upper bound on $l(\setmat)$ immediately translates into algorithms 
for checking primitivity by simply checking every possible product of length smaller or equal to this bound (\rmj{though in some particular cases more efficient algorithms can be used (see Section III)}). 

\subsection{Our results}  

This paper is consequently concerned with upper bounds on the length of $l(\setmat)$ as well as algorithms and complexity of verifying 
existence of a positive product of matrices taken in a given set $\setmat.$ Our main results are: 

\begin{enumerate} \item We show in Section II that recognizing primitivity is decidable but  NP-hard as soon as there are three matrices in the set. Primitivity can be decided in polynomial time for one matrix and so we leave the computational complexity of the case of two matrices unresolved.   

\medskip

\item We also show in Section II that the shortest positive product may have a length that is superpolynomial in the dimension of the matrices, even with a \rmj{fixed} number of matrices in the set.   


\medskip

\item We consider in Section III \aor{the primitivity problem under the additional mild assumption that all the matrices in the set ${\cal P}$ have no zero rows or columns}. We provide a combinatorial proof \rmj{of a \aor{previously-known} primitivity criterion under this assumption.  This resolves an open question of Protasov and Voynov \cite{PV}, who first proved the validity of the \aor{same} criterion \aor{using algebraic tools}}\aor{, and showed it can be
checked in polynomial time.}

\medskip

\item We also prove in Section III that for primitive sets of matrices in ${\cal P}$, the shortest positive product has length $O(n^3)$. Moreover, we show that in this case the length of the shortest positive product is related to the well-known (and unresolved) conjecture of \cerny    on synchronizing automata.  In particular, we show that resolution of the \cerny conjecture would improve the above bound to $O(n^2)$. Moreover, any upper bound on the length of \rmj{the} shortest 
synchronizing word for a {synchronizing} automaton immediately translates into a bound on the length of the shortest positive product of a  
set of primitive matrices in ${\cal P}$. 

\end{enumerate} 

\subsection{Implications of our results} 

\aor{Our results have implications for a number of ongoing research efforts within the field of discrete-time switched systems. First, they complement previous results from \cite{jungers-protasov-lyapunov_ecc,jungers-protasov-lyapunov,protasov-invariant-english,protasov-primitive} which provided simple algorithms for the computation of Lyapunov exponents of nonnegative matrices from ${\cal P}$ for which a positive product exists. If the existence of a positive product is not guaranteed, then the above papers provided more complex and computationally involved protocols relying on quasiconcave maximization. Our results here provide an efficient way of verifying when it is possible to use the lower complexity protocols to compute Lyapunov exponents of matrices from ${\cal P}$. }

\aor{Secondly, our results shed light on the problem of consensus with randomly chosen matrices at each step. Our results in Section III give a necessary and sufficient conditions for primitivity of stochastic matrices (corresponding to consensus on a value in the strict convex hull of the initial states)  which have no zero rows and columns. To our knowledge, the only previous case when necessary and sufficient conditions for consensus with randomly chosen matrices have been provided has been in \cite{jad-r1} for the case of matrices with positive diagonals.  Since stochastic matrices cannot have a zero row by definition, our results  in Section III effectively require only the absence of zero columns, significantly expanding the set of stochastic matrices for which necessary and sufficient conditions for random consensus can be given. } 

\aor{Finally, as we previously remarked, matrix primitivity is among the mathematically simplest possible reachability questions one can pose for switched systems. The NP-hardness results of Section I show that, unfortunately, even this problem cannot be decided in polynomial time unless $P=NP$. In particular, this implies that any generalization of this simple reachability problem is NP-hard as well.} 

\aor{For example, the problem of steering an unknown initial state of Eq. (\ref{eq:switched.system}) to the interior of a given polyhedron by picking the appropriate matrix $A_{i_k}$ at each step is NP-hard, even if the initial condition lies on the boundary of the polyhedron.  More broadly, our results suggest that  reachability problems for discrete-time switched systems can be NP-hard even after a slew of simplifying assumptions on the matrices involved, the structure of the set to be reached, and the initial condition.}

\subsection{Related work} 

The concept of primitive matrix families \rmj{as we study it here} was pioneered in the recent paper \cite{PV}, which extended the classical Perron-Frobenius 
theory and provided a structure theorem  for the primitive matrix sets in ${\cal P}$. A consequence  
of this  theorem was that for matrices in ${\cal P}$ primitivity can be tested in polynomial time. The proofs were based on a somewhat involved spectral 
analysis, and the question of finding a combinatorial proof was left open.  

\rmj{Other generalizations of the well-studied primitivity of one single matrix to a set of matrices exist in the literature.  See for instance \cite{olesky} or \cite{cohen}\aor{, and} \cite{prot-simax} for a recent paper on so-called $k$-primitivity.} 

\ao{We note that two recent papers, appearing simultaneously in 2013 with the conference version of this paper \cite{us-conf}, have also} \rmj{tackled items (3) and (4) above respectively, though with different approaches}\ao{.  
The paper \cite{voynov} proved an $O(n^3)$ bound  
on the length of the  shortest positive product of a set of primitive matrices in ${\cal P}$,  
and  
the paper \cite{alpin} provided a combinatorial proof \aor{of the characterization of primitivity for matrices in ${\cal P}$ first proved in \cite{PV}}. Our work is simultaneous with these results; we remark, however, that our combinatorial proof of item (3) above is shorter relative to \cite{alpin}, and provides \aor{additional} insights.}  \rmj{Finally, our bound in item (4) is stronger than the one proved in \cite{voynov}.}


\section{The general case}  


\rmj{In this section, we study the problem of recognizing primitivity: given a set of $n \times n$ nonnegative matrices $\setmat = \{A_1, \ldots, A_m\}$ does there exist an efficient algorithm that determines whether $\setmat$ is primitive? }

Unfortunately, without any further assumptions on the matrices $A_1, \ldots, A_m$ our main results in  
this section are rather pessimistic. Theorem \ref{hardness} shows that whenever the number of matrices $m$ is at least $3$, testing primitivity is NP-hard; thus there exists no algorithm for recognizing the primivity of a set of $m$ matrices in $R^{n\times n}$ with running time polynomial in $m$ and $n$ unless $P=NP$.  Furthermore, Theorem \ref{lengthseq} shows that the length of the shortest positive 
product can be \aor{superpolynomially} large in the dimension $n$.  

While these results demonstrate that the problem of checking primitivity is intractable in general, we note that 
 it may become tractable under additional assumptions on the matrices $A_1, \ldots, A_m$; indeed,  Section \ref{solvable}  
is dedicated to the study of a class of matrices for which recognizing primitivity has polynomial complexity 
and the length of the shortest possible product is polynomial in $m$ and $n$.  

We now begin with a sequence of definitions and lemmas which will ultimately result in a proof of the aforementioned results, namely Theorem \ref{hardness} on NP-hardness and Theorem \ref{lengthseq} 
on the length of the shortest positive product. We will find it more convenient to make our arguments in terms of graphs rather than matrices; 
our starting point is the following definition which gives a natural way to associate matrices with directed graphs.

\vskip -0.3pc \begin{deff}Given a (directed) graph $G=(V,E)$, the adjacency matrix of $G$, denoted by $A(G)$, is defined as  
\begin{equation*} 
[A(G)]_{ij} = 
\begin{cases} 1, & \text{ if } j \in N_i(G), 
\\ 
0, &\text{ otherwise.} 
\end{cases} 
\end{equation*} where $N_i(G)$ is the out-neighborhood of node $i$ in $G$ (i.e., $N_i(G) = \{ j ~|~ (i,j) \in E\}$).  Conversely, given a nonnegative  matrix $M \in R^{n \times n}$, we will use $G(M)$ to denote the (directed) graph  
with vertex set $\{1, \ldots, n\}$ and edge set $\{ (i,j) ~|~ M_{ij} > 0\}$.  
\end{deff}

It is standard observation that entries of the product $A[G_1] A[G_2] \cdots A[G_l]$ are related to the number of paths in  
the graph sequence $G_1, G_2, G_3, \ldots, G_l$. After formally defining the notion of a path in a graph  
sequence next, we state the relationship in a lemma.

\vskip -0.3pc \begin{deff}Let $G_1, G_2,  \ldots, G_l$ be a sequence of graphs all with the same 
vertex set $V$, and let us adopt the notation $E_k$ for the edge set of  
$G_k$.  For vertices $a,b \in V$ we will say that there exists a path from $a$ to $b$ in  
$G_{1}, \ldots, G_l$ if there exists  a sequence of vertices $i_1, \ldots, i_{l+1}$ such that 
\begin{itemize} \item $i_1 = a$ and $i_{l+1}=b$.  
\item For each $k = 1, \ldots, l$ we have that $(i_k, i_{k+1}) \in E_k$.   
\end{itemize}

We will say that a node $b$ is reachable from $a$ in the sequence $G_1, \ldots, G_l$ if there  
exists a path from $a$ to $b$ in that sequence. Given the graph sequence $G_1, \ldots, G_l$ and node $a \in V$, we will use the notation 
$R_a(k)$ to denote the set of reachable vertices in the sequence $G_1, \ldots, G_k$, where $k \leq l$. We will  
adopt the convention that $R_a(0)=\{a\}$ for all $a \in V$. Finally, we will say that $b$ is reachable from $a$ in $l$ steps of $G_1, \ldots, G_p$ if there exists a sequence of length $l$ consisting of the graphs from 
$\{G_1, \ldots, G_p\}$ in which $b$ is reachable from $a$.  
\end{deff}

The following lemma (which we state without proof and which follows straightforwardly from the  
definition of matrix multiplication) states the usual correspondence between entries of $A[G_1] 
\cdots A[G_l]$ and paths in the sequence $G_1, \ldots, G_l$.

\begin{lem} \label{pathsentries} The number  
of paths from $i$ to $j$ in the sequence $G_1, G_2, \ldots, G_l$  is the $i,j$'th entry of the product $A[G_1] A[G_2] \cdots A[G_l]$. Moreover, there exists a path from $i$ to $j$ in the sequence $G_1, \ldots, G_l$ if and only if 
the $i,j$'th entry of the product  
$P_1 \cdots P_l$ is positive, where for $i=1, \ldots, l$, $P_i$ is any nonnegative matrix satisfying $G(P_i) = G_i$.   
\end{lem}

Thus the primitivity problem for the matrix set $\{A[G_1], \ldots, A[G_k]\}$ is equivalent to the 
problem of finding a sequence of the graphs $G_1, \ldots, G_k$ such that there is at least one path from every node to  
every other node. We will make use of this interpretation shortly.  

We now define the $3$-SAT satisfiability problem, which is well known to be NP-hard.  We will prove below that the primitivity problem is NP-hard by reducing the $3$-SAT problem to it.

\vskip -0.3pc \begin{deff}Let $x_1, \ldots, x_n$ be Boolean variables.  Both   
$x_i$ and its negation $\overline{x_i}$ are called literals. A clause is a disjunction (logical OR)  
of three literals, for example $x_1 \lor x_2 \lor x_3$ or $\overline{x_1} \lor x_2 \lor x_3$. A $3$-CNF formula 
is the conjunction (logical AND) of clauses. For example, $$f =  (x_1 \lor x_2 \lor x_3) \land (\overline{x}_1 \lor \overline{x}_2 \lor \overline{x}_3)$$ is a 3-CNF formula, being a conjunction of 
two clauses. The number of clauses in a $3$-CNF formula is usually denoted by $K$.  Given a 3-CNF formula $f$, the $3$-SAT  
problem asks for an assignment of the values $\{0,1\}$ to the variables $x_1, \ldots, x_n$ so that the 
formula evaluates to $1$. If such an assignment exists, the formula is called satisfiable and the  
corresponding assignment 
is called a satisfying assignment. 
\rmj{The size of an instance of the $3$-SAT problem is $O(n+K).$}
\end{deff} 
Our first step is to associate several graphs with a given $3$-CNF formula,  
as explained in the following definition. \vskip -0.3pc \begin{deff}Given a $3$-CNF formula $f$ on $n$ variables with $K$ clauses, we define three graphs $G_1(f), G_2(f), G_3(f)$.  
Figures \ref{g12} and \ref{g3} show the graphs for the 
the formula $(x_1 \lor x_2 \lor x_3) \land (\overline{x}_1 \lor \overline{x}_2 \lor \overline{x}_3)$.  
We recommend the reader to refer to the figures while going through our  description below.  

All three graphs will have the same vertex set.  
We will have a ``source node'' $u$. For each $i=1, \ldots, K$, we will have the $n$ nodes $u_1^i, \ldots, u_n^i$ and the $n-1$ 
nodes $l_2^i, \ldots, l_n^i$. We will also have the ``failure node'' $f^i$ and the ``success node'' $s^i$; these nodes will also  
be referred to as $u_{n+1}^i$ and $l_{n+1}^i$, respectively.   

For each $i=1, \ldots, K$ and $j = 1, \ldots, n$, if clause $i$ is satisfied by setting $x_j=1$, we put an edge going from $u_j^i$ to 
$l_{j+1}^i$ in $G_1(f)$. Else, we put an edge going from $u_j^i$ to $u_{j+1}^i$ in $G_1(f)$.  

Similarly, for each $i=1, \ldots, K$ and $j = 1, \ldots, n$, if clause $i$ is satisfied by setting $x_j=0$, we put an edge going from $u_j^i$ to 
$l_{j+1}^i$ in $G_2(f)$. Else\footnote{Note that clause $i$ may not contain $x_j$ or its negation $\overline{x_j}$; in that case, neither setting the variable $x_j$ to zero nor to one will satisfy the clause, and  
consequently we will not have a link from $u_j^i$ to $l_{j+1}^i$ in either of $G_1(f)$ and $G_2(f)$.}, we put an edge going from $u_j^i$ to $u_{j+1}^i$ in $G_2(f)$.  

We then add the following edges to both $G_1(f)$ and $G_2(f)$: edges from $l_j^i$ to $l_{j+1}^i$ for all $j=2 \ldots, n$ and $i=1, \ldots, K$; self-loops at all nodes $f^i$; edges leading from each $s^i$ to each $f^i$; and edges leading from $u$ to all $u_1^i$.  

Finally, $G_3(f)$ has edges leading from each $f^i$ and each $s^i$ to $u$, as well as edges leading from each $s^i$ to every node which bears the 
superscript $i$. Note that $G_3(f)$ does not depend on $f$ in the sense that it is the same for all formulas with the same number of variables and same number of clauses.  
\end{deff} 

\begin{center} 
\begin{figure}
\begin{tikzpicture}  
 [scale=.27,auto=left,every node/.style={circle,fill=blue!10},->, thick]
 \node (s)  at (3,10) {$u$}; 
 \node (la)  at (8,4) {$u_1^2$}; 
 \node (lb)  at (8,16) {$u_1^1$}; 
 \node (cb1) at (16,16) {$u_2^1$}; 
 \node (cb2) at (16,12) {$l_2^1$}; 
 \node (rb1)  at (24,16) {$u_3^1$};  
 \node (rb2)  at (24,12) {$l_3^1$}; 
 \node (fb)  at (32,16) {$f^1$};  
 \node (sb)  at (32,12) {$s^1$}; 
 \node (ca1) at (16,4) {$u_2^2$}; 
  \node (ca2) at (16,0) {$l_2^2$}; 
 \node (ra1)  at (24,4) {$u_3^2$};  
 \node (ra2)  at (24,0) {$l_3^2$}; 
 \node (fa)  at (32,4) {$f^2$};  
 \node (sa)  at (32,0) {$s^2$}; 
  
 \path[color=red] 
    (fa) edge [loop above] (fa) 
        (fb) edge [loop above] (fb) 
 (sa) edge [bend left] (fa) 
 (sb) edge [bend left] (fb) 
 (lb) edge (cb1)  
 (rb2) edge [bend right=20] (sb) 
    (cb2) edge [bend right=20] (rb2) 
     (rb1) edge (fb) 
     (cb1) edge (rb1)  
     (la) edge (ca1) 
 (ra2) edge [bend left=20] (sa) 
 (ca2) edge [bend left=20] (ra2) 
 (ra1) edge (sa) 
 (ca1) edge (ra2) 
 (la) edge (ca2)   
 (s) edge [bend left] (la) 
 (s) edge [bend left] (lb); 
  
   \path[color=blue] 
  (fa) edge [loop right] (fa) 
   (fb) edge [loop right] (fb) 
 (sa) edge [bend right] (fa) 
 (sb) edge [bend right] (fb) 
    (rb2) edge [bend left=20] (sb) 
 (cb2) edge [bend left=20] (rb2) 
 (rb1) edge (sb) 
 (cb1) edge (rb2) 
 (lb) edge (cb2)   
   (ra2) edge [bend right=20] (sa) 
    (ca2) edge [bend right=20] (ra2) 
     (ra1) edge (fa) 
     (ca1) edge (ra1)  
     (la) edge (ca1) 
     (s) edge [bend right] (lb) 
     (s) edge [bend right] (la);  
\end{tikzpicture}  
\caption{The graphs $G_1(f)$ (blue) and $G_2(f)$ (red) for the formula $f = (x_1 \lor x_2 \lor x_3) \land (\overline{x}_1 \lor \overline{x}_2 \lor \overline{x}_3)$.} 
\label{g12} \end{figure} 
\end{center} \begin{center}  
\begin{figure} 
\begin{tikzpicture}  
 [scale=.3,auto=left,every node/.style={circle,fill=blue!10},->, thick]
 \node (s)  at (4,10) {$u$}; 
 \node (la)  at (8,4) {}; 
 \node (lb)  at (8,16) {}; 
 \node (cb1) at (16,16) {}; 
 \node (cb2) at (16,13) {}; 
 \node (rb1)  at (24,16) {};  
 \node (rb2)  at (24,13) {}; 
 \node (fb)  at (32,16) {$f^1$};  
 \node (sb)  at (32,11) {$s^1$}; 
 \node (ca1) at (16,4) {}; 
 \node (ca2) at (16,1) {}; 
  \node (ca2) at (16,1) {}; 
 \node (ra1)  at (24,4) {};  
 \node (ra2)  at (24,1) {}; 
 \node (fa)  at (32,4) {$f^2$};  
 \node (sa)  at (32,-1) {$s^2$}; 
  
 \path[color=teal] 
  (fb) edge [bend right=70] (s) 
  (sb) edge [loop below] (sb) 
  (sb) edge (lb) 
  (sb) edge (rb1) 
  (sb) edge (rb2) 
  (sb) edge (cb1) 
  (sb) edge [bend left](cb2) 
  (sb) edge (fb) 
  (sb) edge (s) 
  (fa) edge [bend left=70](s) 
  (sa) edge [loop below] (sa) 
  (sa) edge (la) 
  (sa) edge (ra1) 
  (sa) edge (ra2) 
  (sa) edge (ca1) 
  (sa) edge [bend left](ca2) 
  (sa) edge (fa) 
  (sa) edge (s); 

\end{tikzpicture}    
\caption{The graph $G_3(f)$ for the formula $f = (x_1 \lor x_2 \lor x_3) \land (\overline{x}_1 \lor \overline{x}_2 \lor \overline{x}_3)$.} 
\label{g3}  \end{figure}
\end{center} \vskip -2.8pc We remark that this construction is a variation of one of the constructions from the earlier work \cite{tsitsiklis97lyapunov}.  It appears somewhat unwieldy at first glance, but the subsequent lemmas will provide 
some insights into it. First, however, we state our first main result of this section, Theorem \ref{hardness}, which provides a reduction from $3$-SAT to checking primitivity of  
a set of three matrices.

\begin{thm} \label{hardness} The $3$-SAT formula $f$ has a satisfying assignment if and only if  
the matrix set $$\{A[G_1(f)], A[G_2(f)], A[G_3(f)]\}$$ is primitive. Consequently, there is no  
algorithm for deciding matrix primitivity which scales polynomially in $n$ unless $P=NP$.  
\end{thm}

We now begin the proof of this theorem. We will assume henceforth that $f$ is a fixed formula, and  
correspondingly we will simply write $G_1, G_2, G_3$ for the three graphs. We begin with the key 
lemma which encapsulates the most important property of these graphs. We remark that this is a variation 
of a lemma from \cite{tsitsiklis97lyapunov,bl-ol} used to establish the complexity of closely related problems.

\begin{lem} \label{sat-paths} Consider a sequence of length $n$ of graphs from $\{G_1,G_2\}$ and  
set $x_i=1$ if the $i$'th 
graph is $G_1$, and $x_i=0$ if the $i$'th graph is $G_2$. We have that $s^i$ is reachable from 
$u_1^i$ in this sequence if and only if the $i$'th clause of $f$ is satisfied by the 
assignment $x_1, \ldots, x_n$.  
\end{lem} 

\vskip -1.5pc \begin{pf} By construction, an edge goes from $u_j^i$ to $l_{j+1}^i$ in $G_1$ whenever setting $x_j=1$ satisfies the $i$'th clause of  
$f$, and the same edge is present in $G_2$ whenever setting $x_j=0$ satisfies the clause of $f$. Thus  
clause $i$ is satisfied by the assignment of $x_1, \ldots, x_n$ defined in this  
lemma if and only if the corresponding sequence of $G_1,G_2$ includes an edge from some $u_j^i$ to $l_{j+1}^i$. But then the presence of edges from each 
$l_j^i$ to $l_{j+1}^i$ implies this happens if and only if $s^i = l_{n+1}^i$ is reachable 
after $n$ steps from $u_1^i$.   \qedsymbol
\end{pf}

This simple lemma is an important ingredient of our proof of Theorem \ref{hardness}. Indeed,  
to prove this theorem we need to relate the satisfiability of $f$ to the primitivity of the 
matrix set $\{A[G_1], A[G_2], A[G_3]\}$. The latter, as a consequence of Lemma \ref{pathsentries}, 
can be recast as a question about the existence of a sequence with a path from every node to every other node; we thus need to somehow 
relate path-existence questions to satisfiability questions. This is precisely what is done by this 
previous lemma.  

We sharpen the conclusions of this lemma with the following corollary, which follows from the  
fact that $\{u_1^1, \ldots, u_1^K\}$ is the set of out-neighbours of $u$:

\begin{cor} Consider a sequence of length $n+1$ of graphs from $\{G_1,G_2\}$, and define $x_i=1$ if the $i+1$'st 
graph is $G_1$, and $x_i=0$ if the $i+1$st graph is $G_2$. We have that all $s^i$ are reachable from 
$u$ in this sequence if and only this $x_1, \ldots, x_n$ is a satisfying assignment for $f$. \label{upaths} 
\end{cor}

We are now essentially ready to provide a proof of Theorem \ref{hardness}. However, before embarking  
on the details of the proof, we collect a number of straightforward observations about the  
graphs $G_1,G_2, G_3$ in a remark.

\begin{rem} \label{obviousremark} $~$

\begin{itemize}  
%
%
%
%
%

\item The set of reachable nodes from $u$ in strictly more than  
$n+1$ steps of $G_1,G_2$ is a nonempty subset of the  failure nodes $f^i$.

Indeed, this follows from the previous item and the observation that the only outgoing link from $s^i$ in these graphs leads to 
$f^i$, and the only outgoing link from $f^i$ leads to itself.

\item The nodes reachable from a node $v \neq u$ in $n+1$ steps of $G_1,G_2$ or more are a nonempty subset of the set 
of failure nodes.

The argument for this is identical to the argument for the previous item.

\item Consider the sequence $G_1, G_1, \ldots$ repeated $n+1$ times. Regardless of the 
starting vertex, the only reachable vertices are success nodes and failure nodes.

This follows by the previous item and the fact that $G_3$ has a single outgoing edge from each  
failure node to $u$.  
\end{itemize} 
\end{rem}

\vskip -1.5pc \begin{pf}[Proof of Theorem \ref{hardness}] Suppose first that the $3$-SAT problem has a satisfying assumption. Consider the graph sequence 
of length  $2n+4$ defined as follows:

\begin{enumerate} \item  First, we repeat the graph $G_1$ $n+1$ times. 
\item  The $n+2$'nd graph 
equals $G_3$ and the $n+3$'rd graph is $G_1$  
\item For $k=1, \ldots,n$ the $n+3+k$'th graph is $G_1$ if $x_k=1$ in the satisfying assignment and $G_2$ if 
$x_k=0$ in the satisfying assignment.  
\item Finally, $G_{2n+4}=G_3$.  
\end{enumerate}

For example, for the formula   
$(x_1 \lor x_2 \lor x_3) \land (\overline{x}_1 \lor \overline{x}_2 \lor \overline{x}_3)$, we have the 
satisfying assumption $x_1=0, x_2=0, x_3=1$; the corresponding sequence is $$G_1, G_1, G_1, G_1, G_3, G_1,  
G_2, G_2, G_1, G_3.$$  

We claim that there are paths from every node to every other node in this sequence, which by Lemma  
\ref{pathsentries} implies that the corresponding matrix product is positive. Indeed, as noted  
in Remark \ref{obviousremark}, after repeating $G_1$ $n+1$ times, we have that for any node $a \in V$, $R_a(n+1)$ is a nonempty subset of the set of failure nodes and success nodes. After graph $G_3$ on the $n+2$'nd step, we 
have $u\in R_a(n+2);$ and after $G_1$ on the $n+3$'rd step, we have  
$R_a(n+3) = \{u_1^1, u_1^2, \ldots, u_K^1\}$.  

Now appealing to Lemma \ref{sat-paths}, we see that $R_a(2n+3)$ contains the set of all 
the success node $\{s^1, \ldots, s^n\}$. When we apply the  
graph $G_3$ on the $2n+4$'th step, it follows that every node becomes reachable.

Conversely, suppose that there exists a positive product of $A[G_1], A[G_2], A[G_3]$. By Lemma   
\ref{pathsentries}, there exists a sequence of $G_{1},G_{2}, G_{3}$ such that, in particular, 
there is a path from $u$ to every other node. Consider one such sequence of minimal length; say it has 
length $p$.  

Note that since neither $G_1$ nor $G_2$ have an edge incoming 
to $u$, it follows that the $p$'th graph in this sequence must be $G_3$. Consider the last time $G_3$ appeared {\em before} time $p$. Let us say that this happened at time $k$, i.e., the $k$'th graph was $G_3$ and $G_3$ did not appear in the sequence at any time $l$ which satisfies $k < l < p$. In the event that $G_3$ never 
appeared before time $p$, we will set $k=0$.  

Now starting from $u$ the set of reachable nodes $R_u(k)$ is certainly non-empty: this is trivially true 
if $k=0$, and otherwise true because $R_u(p)=V$ by assumption. Consequently, due to the structure of $G_3$, we have $u \in R_u(k)$. Moreover, since $R_u(k) \neq V$ (else $p$ would not be the shortest length of a positive 
product), the structure of $G_3$ implies that some $s^i$ does not belong to $R_u(k)$. Without loss of  
generality, suppose $s^1 \notin R_u(k)$. Yet once 
again appealing to the structure of $G_3$,  
this implies all nodes with superscript $1$ do not belong to $R_u(k)$.  

We next argue that $p-1-k=n+1$, i.e., there are exactly $n+1$ graphs between $k$ and $p-1$. Indeed, we know that $R_u(p)=V$; but in order for $s^1 \in R_u(p)$, by the structure of $G_3$ we must 
have $s^1 \in R_u(p-1)$. Since in any sequence of $G_1,G_2$ (i)  the only node without a superscript of $1$  
that has a path to $s^1$ of any length is $u$ (ii) the only path from $u$ to $s^1$ has length $n+1$, we must have that $(p-1) - k = n+1$, as claimed.  

Furthermore - once again by the structure of $G_3$ - we must have that every $s^i$ belongs to  
$R_u(p-1)$ in order to have $R_u(p)=V$. As  pointed out earlier in Remark \ref{obviousremark}, starting 
from any node other than $u$ there are no paths of length $n+1$ to any $s^i$. It follows that in the  
minimal length sequence we are considering there 
must be a path of length $n+1$ from $u$ to {\bf all} $s^i$.  

Now we appeal to Corollary \ref{upaths} to get a satisfying assignment for $f$.  \qedsymbol
\end{pf}

We now turn to the question of bounding $l(\setmat)$, the length of the shortest positive product of matrices from $\setmat$; we will adopt the convention that $l(\setmat) = +\infty$ when no product of matrices from  
$\setmat$ is positive. As we remarked earlier, any upper bound on $l(\setmat)$ can be translated into an algorithm for checking primitivity simply by checking all products of length $l(\setmat)$ from $\setmat$. Unfortunately, our results are once again quite pessimistic: while upper bounds exist that show matrix 
primitivity is decidable, we construct four nonnegative matrices for  
which the shortest positive product has length at least exponential in the dimension.  

We define $l(m,n)$ to be the largest $l(\setmat)$ over all sets $\setmat$ with $m$ matrices of size $n \times n$ with $l(\setmat)<\infty$. Our second main result of this section is the following theorem.

\begin{thm}  \label{lengthseq} We have that for all $m,n$,  \[ l(m,n) \leq 2^{n^2}.\] Moreover, if $m \geq 4$, then for all $\epsilon>0$ there exists a sequence of positive integers $n_1, n_2, \ldots,$ tending to infinity such that 
\[ \left((1-\epsilon) e \right)^{\sqrt{n_k/2}} \leq l(m,n_k). \] 
\end{thm}

We see that $l(m,n)$ does not in general strongly depend on the number of matrices $m$, in the sense that in 
the case of $m \geq 4$ it can be bounded 
above and below by exponentials independent of $m$. An obvious consequence of this theorem is that  
matrix primitivity is decidable, but the  natural algorithm which tests all products of  
length $l(m,n)$ can take doubly exponential time in the dimension to halt.  

We conclude this section with a proof of this theorem.

\vskip -1.5pc \begin{pf}[Proof of Theorem \ref{lengthseq}] We first do the easy direction, namely we prove  
inequality $$l(m,n) \leq 2^{n^2}.$$ Let $P_1 P_2 \ldots P_k$ be a  
shortest-length positive product of matrices from $\setmat$. For $S \subset \{1, \ldots, k\}$,  
we will use the notation $P_S$ to denote the product  
\[ P_S = \prod_{i \in S} P_i,\] where the indices are taken in increasing order. We argue that we cannot have  
\begin{equation} \label{requal} R_1(i_1) = R_1(i_2), ~~ R_2(i_1) = R_2(i_2), ~~ \ldots, ~~ R_n(i_1) = R_n(i_2) \end{equation} for some $1 \leq i_1 < i_2 \leq k$. Indeed, let us proceed by contradiction; 
assume Eq. (\ref{requal}) holds for such $i_1,i_2$;  we then argue that the product $P_{\{1,\ldots, k\} \setminus \{i_1+1, i_1+2, \ldots, i_2 \}}$ is also positive. Indeed, we know  
by Lemma \ref{pathsentries} that positivity 
of the product $A_1, A_2, \ldots, A_l$ is equivalent to the existence of a path from every $i$ to every $j$ in the sequence 
$G(A_1), G(A_2), \ldots, G(A_l)$. But since Eq. (\ref{requal}) implies that the set of reachable points by times $i_1$ and $i_2$ is the  
same, this means there is a path from every $i$ to every $j$ in the subsequence of $G(P_1), \ldots, G(P_k)$ which omits the graphs $G_{P_{i_1+1}}, \ldots, G_{P_{i_2}}$.  
Consequently, $P_{\{1,\ldots, k\} \setminus \{i_1+1, i_1+2, \ldots, i_2 \}}$ is positive. However,  $P_1 \cdots P_k$ was assumed to be a minimum length positive product; we therefore conclude that indeed Eq. (\ref{requal}) cannot hold.  

This means that the $n$uple $(R_1(i), R_2(i), \ldots, R_n(i))$ takes on distinct values for each $i=1, \ldots, k$.  But since each set $R_i$ can assume at most $2^n$ possible values, the tuple has at most $2^{n^2}$ possible 
values and thus $k$ cannot be larger than $2^{n^2}$. This proves the upper bound.

We now turn our attention to the lower bound of the theorem. Let $n_k = 1+\sum_{i=2}^k p(i)$ where $p(i)$ is the $i$'th prime and $k \geq 3$. We will establish the lower bound as follows: we will describe four graphs  
on $n_k$ nodes whose adjacency matrices have a positive product, but the shortest length positive product has length at least $\left( (1-\epsilon) e \right)^{\sqrt{n_k/2}}$ as long as $n_k$ is large enough as a  
function of $\epsilon$.  

We next describe the graphs.  We refer the reader to Figures \ref{g12length} and \ref{g3length}, and we recommend that the reader refer to  figures in while following our description.  We will have a source node $u$ and $k-1$ cycles, the $i$'th on $p(i)$ nodes. $G_1$ will have edges going around each cycle counterclockwise. $G_2$ will have edges from the last node in each cycle to all the nodes in its cycle, as well as an edge from the last node in each cycle to $u$. $G_3$ will have edges from every node to $u$, and $G_4$ will have edges from $u$ to the first node in each cycle. We let $\setmat$ be $\{A[G_1], A[G_2], A[G_3], A[G_4]\}$, i.e., the adjacency matrices of these four graphs.

We first show that $l(m,n_k)$ is finite by exhibiting a positive product. We pick the following sequence of graphs. The first is $G_3$, which ensures that 
$R_i(1)=u$ for any $i$; the second is $G_4$ which ensures that $R_i(2)$ is exactly the set of the first nodes in all cycles for any $i$. We then pick $G_1$ and repeat it $P(k) = 2 \cdot 3 \cdot 5 \cdots p(k)-1$ times.  Observe that (i) for each $i=2, \ldots, k$, $P(k) - (p(i)-1)$ is divisible by $p(i);$ and (ii) starting from the first node in each cycle, $p(i)-1$ steps bring us to the last node in that cycle, and then 
a multiple of $p(i)$ steps returns us to the same node in that cycle. Therefore, the set of nodes reachable (from any node)  
after the last repetition of $G_1$ is exactly the set of the last nodes in each cycle. Applying $G_2$ then means that the set of reachable nodes equals the set of all nodes. We have just exhibited a sequence 
with the property that any node is reachable from any node; by Lemma \ref{pathsentries} the corresponding matrix product of $A[G_1], A[G_2], A[G_3], A[G_4]$ is positive.

We now proceed to prove a lower bound on the length of the shortest positive product of $A[G_1], A[G_2], A[G_3], A[G_4]$. Consider such a product; it corresponds to a sequence of graphs from $\{G_1, G_2, G_3, G_4\}$ of shortest length in which there is a path from 
any node to any other node. As we have argued previously, the fact that this sequence has minimal length  
implies the tuple of reachable sets never repeats.  

We can immediately assert that the first graph is $G_3$; if any other graph 
appears as first then some $R_i(1)$ equals the empty set, and so equals the empty set thereafter.  
Similarly, the second graph then cannot be $G_1$ or  
$G_2$ because that would give $R_i(2)=\emptyset$ for all $i$; it therefore must either be $G_3$ or $G_4$. 
It cannot be $G_3$ because that would contradict the minimiality of the sequence; so it is $G_4$.  

We thus have that $R_i(2)$ is the set of the first nodes of every cycle for every node $i$. We now argue that $G_3$ never occurs again, since that would make the reachable set from any node $i$ equal to $R_i(1)$ and 
once again contradict the minimality of the graph sequence we are considering. Thus the remainder of the sequence is composed of just 
$G_1,G_2,G_4$.  

We now argue that in the remainder of the sequence $G_2$ appears only once, and on the last step.  
Indeed, it is clearly true that $G_2$ must appear on the last step since each of the other graphs has at 
least one node without any incoming edges. On the other hand, suppose $G_2$ appears before the last  
step, at time $k$. Now either we have that for all $i$, $R_i(k-1)$ includes all of the last nodes in all the cycles, or some  
$R_i(k-1)$ does not include some last node in some cycle.  
In the former case, $R_i(k)$ is the set of all nodes for every $i$, which contradicts the  
minimality of the sequence. In the latter case, there is at least one cycle such that all $R_i(k)$ do not 
include any node in that cycle; without loss of generality, suppose it is the cycle of length $2$.  This means the graph $G_4$ appears at some future time as it is the only graph with an edge incoming to a node on the length 
$2$ cycle. But after applying $G_4$ the set of reachable nodes is either empty or consists of all the first nodes  
in all the cycles, both of which cannot be: the first contradicts the eventual existence of paths from every node 
to every node, and the second contradicts the minimality of the sequence since this tuple of reachable 
sets has already occured.  

We recap: we have shown that the sequence is of the form $G_3, G_4, *, \ldots, *, G_2$ where $*$'s are  
either $G_1$ or $G_4$. But we can now immediately conclude that every $*$ above is in fact $G_1$: if some $*$ were $G_4$, the set of 
reachable nodes is empty.  

We have thus concluded that the sequence must be of the form $$G_3, G_4, G_1, \ldots, G_1, G_2.$$ To be  
able to bound its length, we have to consider just how many repetitions of $G_1$ this sequence has.  Note that $R_i(2)$ consists of the first 
nodes in each cycle, while after the last $G_1$, we have that $R_i$ has to include the last node in each cycle.  
This means that the number of times $G_1$ is repeated has to be of the form $p(i)-1 + j(i) p(i)$ for each $i=2, \ldots, k$, where $j(i)$ 
is a nonnegative integer: $p(i)-1$ times to reach the last node in the cycle and $j(i) \geq 0$ additional 
round trips of length $p(i)$ around the cycle.  

Note that since $$p(2)-1+j(2) p(2) = p(3) - 1 + j(3) p(3)$$ we have  
$$ (1 + j(2)) p(2) = (1 + j(3)) p(3)$$ and consequently $p(3)$ divides $1+j(2)$. Repeating this 
argument, we see that $$ 1 + j(2) \geq \prod_{l=3}^{k} p(i)$$ so that the number of times $G_1$ is 
repeated is at least $$p(2) - 1 + \left( -1 + \prod_{l=3}^k p(i) \right) p(2) = -1 + \prod_{l=2}^k p(i).$$ 
Thus the total length of the sequence is at least  
$2 + \prod_{l=2}^k p(i)$ while the number of nodes is $1 + \sum_{i=2}^k p(i)$. Now we can conclude using the  
prime product limit \cite{ruiz},   
\[ \lim_{k \rightarrow \infty} \left( \prod_{l=2}^k p(i) \right)^{1/p(k)} = e \] and the fact that for 
any $\epsilon>0$, we have that for large enough $i$, \[ i \ln i \leq p(i) \leq (1+\epsilon) i \ln i \] which implies that for any $\epsilon>0$ for large enough $k$ \cite{bach-shalitt}, \[ \prod_{l=2}^k p(i) \geq \left( (1-\epsilon) e \right)^{k \log k} \] and also that (again for large enough $k$), 
\[ 1 + \sum_{i=2}^k p(i) \leq 2 k^2 \log k . \] Thus with $n_k= 1 + \sum_{i=2}^k p(i)$ we have  
shown that the minimal length of a positive product is finite and at least $\left( (1-\epsilon) e \right)^{\sqrt{n_k/2}}$ as long as $n_k$ is large enough as a  
function of $\epsilon$; this concludes the proof of this theorem. \qedsymbol
\end{pf} 
 
\begin{center}  
\begin{figure} 
\begin{tikzpicture}  
 [scale=.36,auto=left,every node/.style={circle,fill=blue!10},->, thick] 
 \node (s)  at (0,10) {$u$}; 
 \node (la)  at (6,4) {}; 
 \node (lb)  at (6,16) {}; 
 \node (cb2) at (14,16) {}; 
 \node (cb3) at (14,4) {}; 
 \node (cb4) at (10, 9) {};

 \path[color=blue] 
    (lb) edge [bend left] (cb2) 
    (cb2) edge [bend left] (lb) 
    (la) edge (cb3) 
    (cb3) edge (cb4) 
    (cb4) edge (la); 

    \end{tikzpicture}  \caption{Proof of Theorem \ref{lengthseq}: the graph $G_1$} 
    \label{g12length} \end{figure}
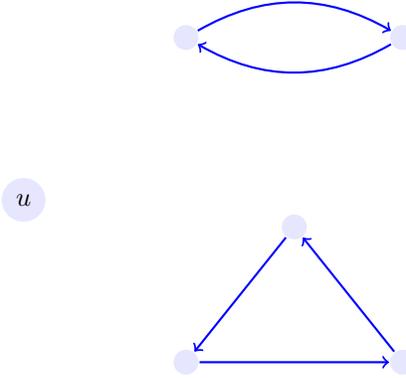 
\end{center} 

\begin{center}  
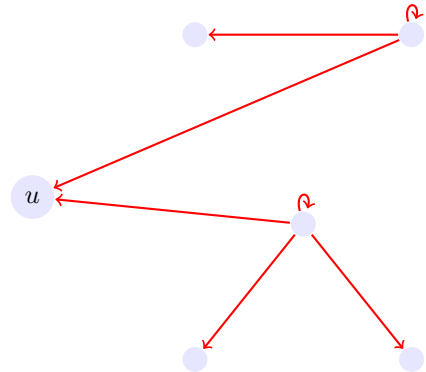
\begin{figure} 
\begin{tikzpicture}  
 [scale=.36,auto=left,every node/.style={circle,fill=blue!10},->, thick]

   \node (s2)  at (20,10) {$u$}; 
 \node (la2)  at (26,4) {}; 
 \node (lb2)  at (26,16) {}; 
 \node (cb22) at (34,16) {}; 
 \node (cb32) at (34,4) {}; 
 \node (cb42) at (30, 9) {};

  \path[color=red] 
   (cb22) edge (lb2) 
   (cb22) edge [loop above] (cb22) 
   (cb42) edge [loop above] (cb42) 
   (cb42) edge (cb32) 
   (cb42) edge (la2) 
    (cb22) edge (s2) 
   (cb42) edge (s2); 
    \end{tikzpicture}  \caption{Proof of Theorem \ref{lengthseq}: the graph $G_2$} 
    \label{g3length} \end{figure} 
\end{center} 

This concludes our section on general nonnegative  matrices.  We have proved that not only it is NP-hard to decide whether a set is primitive, but even more, the size of such a minimal product can be exponential in the size of the matrices.  In the next section we show that under an additional assumption, the situation changes, and primitivity becomes much more amenable on an algorithmic point of view. 

\section{Sets with no zero rows nor zero columns \label{solvable}} 
In this section we focus on sets of matrices that satisfy the following assumption: 
\begin{assumption}\label{assum-nzrnzc} 
No matrix $A \in \setmat$ has a row or a column identically equal to zero: 
$$\forall A \in \setmat,\, \forall 1\leq i \leq n,\, \exists j: A_{i,j}>0, $$ 
$$\forall A \in \setmat,\, \forall 1\leq j \leq n,\, \exists i: A_{i,j}>0. $$ 
\end{assumption} 

We start with an easy lemma. In the following, we write $\setmat^t$ for the set of matrices which are products of length $t$ of matrices taken in $\setmat,$ and $\setmat^*$ for the set of products of arbitrary length of matrices in $\setmat.$ 

\begin{lem}\label{lem-assumption1} 
If a set of matrices $\setmat$ satisfies assumption \ref{assum-nzrnzc}, then every matrix in $\setmat^*$ satisfies it. 
\end{lem} 

\subsection{\rmj{A combinatorial proof for the Protasov-Voynov characterization}} 

The following theorem provides a structural characterization of primitivity for sets of matrices satisfying Assumption \ref{assum-nzrnzc}.  It was first proved in \cite{PV} (after a conjecture of \cite{protasov-primitive}), where the authors show that it leads to an efficient algorithm for recognizing such sets. The proof in that paper is long, and involves linear algebraic and geometric considerations.  The authors also ask whether a simple combinatorial proof is possible for this result.  We provide here an alternative combinatorial and self-contained proof. \rmj{As mentioned in the introduction, see \cite{alpin} for a related work simultaneous and independent to ours.}
\begin{thm}\cite{PV} 
A set of nonnegative matrices $\setmat \subset \re^{n\times n}_+$ satisfying Assumption \ref{assum-nzrnzc} fails to be  primitive if and only if one of the following conditions holds: 
\begin{itemize} 
\item There exists a permutation matrix $P$ such that all the matrices $A\in \setmat$ can be put in the same block triangular structure.  Equivalently, there exists a partition\footnote{A \emph{partition} of a set $S$ is a set ${\cal{ P}}=\{N_1,N_2,\dots,N_k\},\, k>1,$ such that \ao{all $N_i$ are nonempty,} $N_i \subset \mathbb S,\, N_i\bigcap N_j=\emptyset\  \forall i \neq j,\, \mbox{ and } N_1\bigcup N_2\dots \bigcup N_k = S.$ }$\{N_1,N_2\}$ of $\{1,\dots,n\}$ such that 
\begin{eqnarray}\forall A\in\setmat, i\in N_1, j\in N_2 \Rightarrow A_{i,j}=0.& 
\end{eqnarray} 

\item There exists a permutation matrix $P$ such that all the matrices $A\in \setmat$ can be put in the same block permutation structure. Equivalently, there exists a partition $\{N_1,N_2,\dots,N_k\}$ such that for all $ A\in\setmat,$ 

\begin{eqnarray}\nonumber 
\exists i\in N_I,j\in N_J, A_{i,j}\neq 0 \Rightarrow\\ \nonumber \forall i'\in N_I,j'\in N_J,k\not\in N_I, l\not\in N_J,\\ A_{i',l}=0, A_{k,j'}=0. 
\end{eqnarray} 
\end{itemize} 
\end{thm} 
\vskip -1.5pc \begin{pf} 
The two conditions in the theorem are obvious sufficient conditions for imprimitivity.  Indeed, it is easy to see that the product of two matrices satisfying any of these conditions still satisfies it, and thus cannot be primitive.

Suppose now that the set is imprimitive, and the first condition is violated (that is, the set is irreducible).  We prove that the second condition holds.   

Let then $n_1\leq n$ be the maximum number of nonzero entries in a row or a column of any matrix in $\cm^\ao{*}.$

We {\bf claim} that, up to relabelling of the entries, there is a product $A \ao{\in \cm^\ao{*}}$ such that $\forall 1\leq i,j \leq n_1,$ $A_{i,j}>0.$  {This implies that $n_1<n.$}

\ao{Indeed}, let us take a product $B\in \cm^*$ \ao{which has a row or a column with $n_1$ positive entries. Let us suppose that it is a row, and let $k$ be the index of this row. Let us further denote by $S$ the set of indices \ao{coresponding to} these positive entries:$$  B_{k,i}>0 \Longleftrightarrow i \in S.$$ Without loss of generality, } we can futher assume that \rmj{$B_{k,k}>0$} - if not, we can simply pre-multiply $B$ by a matrix $C$ in $\setmat^*$ such that $C_{k',k}>0$ for some $k'\in S$ and this is always possible since $\setmat$ is irreducible. \ao{Observe that $B_{ij} = 0$ for all $i \in S, j \in \rmj{V\setminus S}$: else $B^2$ would have $n_1+1$ nonzero entries in the $k$'th row.}

\ao{Note that $S' = \{k\}$ is a subset of $S$ such that $B_{ij} > 0$ for all $i \in S', j \in S$. We argue that for any such $S'$ which is a strict subset of $S$, we can construct $S''$ with the same property whose cardinality is one larger than the cardinality of $S'$. Indeed, let $k \in S \setminus S'$ and let $r$ be such that $B_{kr}>0$. As we argued above, $r$ must be in $S$. Let $q$ be any element of $S'$ and take $D$ to be a matrix with $D_{rq}>0$ - such a $D$ exists by irreducibility. Then the matrix $BDB$ has the property that $B_{ij} >0$ for all $i \in S' \cup \{k\}, j \in S$, proving the claim.}

\ao{Iterating this argument, we obtain the existence of the $A \in \cm^*$ such that $A_{ij} >0 $ whenever $i,j \in S$, proving the boldfaced claim above. A similar argument works if $B$ has a column, rather than a row, of $n_1$ positive entries. }

\ao{Observe that} since $n_1$ is the maximum number of nonzero entries in any row and column, the matrix $A$ is actually block-diagonal.


We now iterate this argument inductively to obtain a  partition $\{S_1,\dots, S_p\}$ of $\{1,\dots, n\}$. Indeed, define \ao{$S_q$ as the largest set within $V \setminus \left( S_1 \cup \cdots \cup S_{q-1} \right)$ such that there exists a product with entries corresponding to $S_q$ positive in some row or column.} We claim that the following two properties are satisfied:

\begin{enumerate}  
\item  There exists a matrix $A_q\in \cm^*$ such that \ao{ $[A_q]_{ij} > 0$  if $i,j \in S_q$ and $[A_q]_{ij} = 0$ if one of $i,j$ belongs to $S_q$ and the other does not.} 

\smallskip

\item \ao{There does not exist a matrix $B_q \in \cm^*$ with $[B_q]_{ij}>0$} {for all} \ao{$i,j \in S_q$ and   $[B_q]_{kl} > 0$ for some $k,l,$ one of which is in $S_q$ and the other is not. }\end{enumerate}  



\ao{Observe that we have already proven the ``base case'' of $S_1$. Assume now we have proven the case of $S_1, \ldots, S_{q-1}$ and consider $S_q$.} \ao{Let $B'$ be the product with $|S_q|$ nonzero entries in a row or column in the definition of $S_q$. Suppose it is a row in $B'$ with entries corresponding to $S_q$ positive. Let it be the $k$'th row; as before, we may assume by irreducibility that $[B']_{kk}>0$. We then argue that, as before, $[B']_{ij} = 0$ if $i \in S_{q}$ and $j \notin S_{q}$. Indeed, if $[B']_{ij}  > 0$ for $j$ in some $S_{q'}$ with $q' > q$, then $[B']^2$ has more than $n_q$ positive entries in $V \setminus \left( S_1 \cup \cdots \cup S_{q-1} \right)$, a contradiction. Alternatively, if $q' < q$, then taking any $a \in S_{q'}$ and $D$ any matrix with $D_{ak}>0$ we have that $A_{q'} D (B')^2 A_{q'}$} \rmj{contradicts} \ao{item (2) of the inductive hypothesis: it has positive entries for all rows and columns corresponding to $S_{q'}$, as well as at least one more positive entry in each row corresponding to $S_{q'}$.}\\
\ao{We thus have that $S'=\{k\}$ satisfies $[B']_{ij} > 0$ whenever $i \in S', j \in S_q$. This implies the existence of the matrix $A_q \in \cm^*$ which satisfies $[A_q]_{ij} > 0$ whenever $i \in S_q, j \in S_q$ by repeating verbatim the corresponding steps for $S_1$ peformed earlier in the proof.}

\ao{It remains to argue that no matrix $B \in \cm^*$ (including the matrix $A_q$) can have $B_{ij}>0$ for all $i,j \in S_q$ as well as $B_{kl}>0$ when one of $k,l$ belongs to $S_q$ and the other does not. Indeed, suppose for example that $k \in S_q, l \in S_{q'}$. If $q'>q$, the definition of $S_q$ is contradicted. If $q'<q$, then taking any $a \in S_{q'}$ and $D$ such that $D_{ak}>0$ we have that $A_{q'} DB A_{q'}$ violates item (2) of the inductive hypothesis. The case when $k \notin S_q, l \in S_{q'}$ is similar, as is the case when it is a column of $B'$ that has $|S_q|$ positive entries is similar. This concludes the proof of items (1) and (2) above. }

\ao{We can now argue that every matrix in the set $\cm$ is a permutation on the sets $\{S_1,\dots, S_p\}.$  We must prove two statements, one to the effect that links ``out-going'' from the same $S_i$ cannot lead to different $S_i$'s, and one to the effect that links ``in-coming'' to the same $S_i$ cannot come from distinct $S_i$'s. We will prove the former statement (and the proof of the latter is similar). Formally, we argue that it is impossible to have $B \in \cm^*$ and $i,k \in S_{q}$ such that $B_{ij}>0, B_{kl}>0$ with $j \in S_{q'}$ and $l \in S_{q''}$ and $q' \neq q''$.  }

Indeed, suppose such a matrix \ao{$B \in \cm^*$} \ao{exists.}
Then, taking a matrix $C\in \cm^*$ such that $C_{v,w}>0$ for some $v\in q',w\in q,$ \ao{the product $A_{q'}CA_{q}BA_{q'}$ violates item (2) above. Indeed, it has $n_{q'}+1$ positive entries in all the rows of $S_{q'}$, corresponding to all the columns of $S_{q'}$ as well as at least one more. }
\qedsymbol \end{pf}

\subsection{Bounds on the length of the product} 
We now turn to the problem of obtaining tight bounds on the length of a shortest strictly positive product, as a function of the dimension of the matrices.  
For this purpose, we make connections with a well known concept in TCS, namely, Synchronizing Automata.  Our result also suggests that an exact answer to that problem is probably very hard to obtain. 

A (deterministic, finite state, complete) automaton is a set of $m$ row-stochastic matrices $\setmat \subset \{0,1\}^{n\times n}$ (where $m,n$ are positive integers).  That is, the matrices in $\setmat$ have binary entries, and they satisfy  $ A \vectun =\vectun,$ where $\vectun$ is the all-ones (column) vector.   For convenience of product representation, to each matrix $A_c\in \setmat$ is associated a letter $c, $  such that the product $A_{c_1}\dots A_{c_t}\in \setmat^t$ can be written $A_{c_1\dots c_t}.$ 

\vskip -0.3pc \begin{deff}
An automaton $\setmat \subset \{0,1\}^{n\times n}$ is \emph{synchronizing} if there is a finite product $A=A_{c_1}\dots A_{c_T}:\, A_{c_i}\in \setmat$ which satisfies $$A=\vectun e_i^T, $$ 
where $\vectun$ is the all-ones vector and $e_i$ is the $i$th standard basis vector.\\ 
In this case, the sequence of letters $c_1\dots c_T$ is said to be a \emph{synchronizing word.} 
\end{deff} 
We recall the following conjecture which has raised a large interest in the TCS community \cite{kari03}, \cite{eppstein90reset}, \cite{trahtman_cerny},\cite{jungers_sync_12},\cite{cerny64}.  It has been proved to hold in many particular cases, but the general case remains open. 
  
\begin{figure} 
\centering 
\begin{tabular}{c} 
\includegraphics[scale = 0.2]{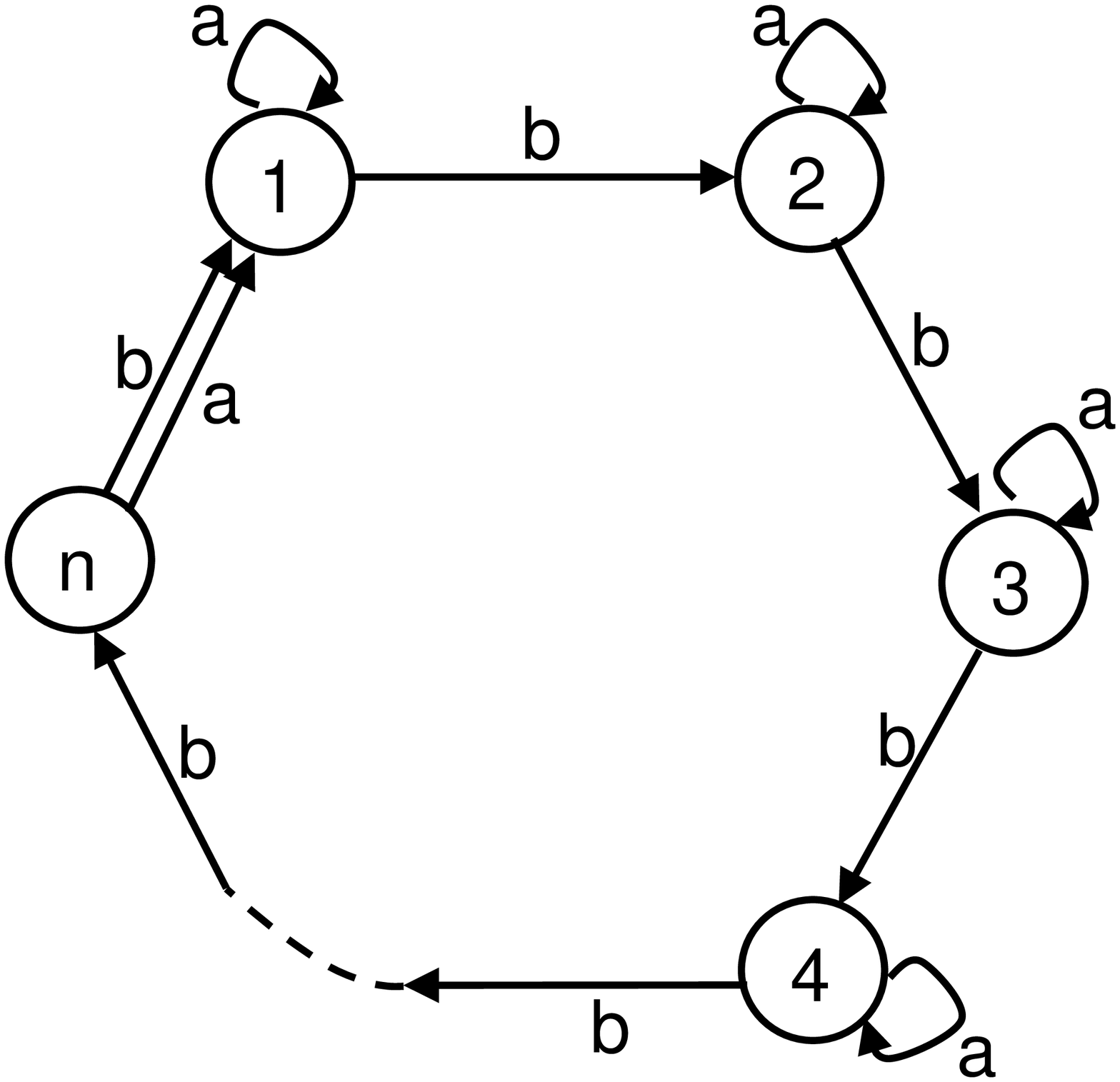}\\(a)\\ 
\includegraphics[scale = 0.2]{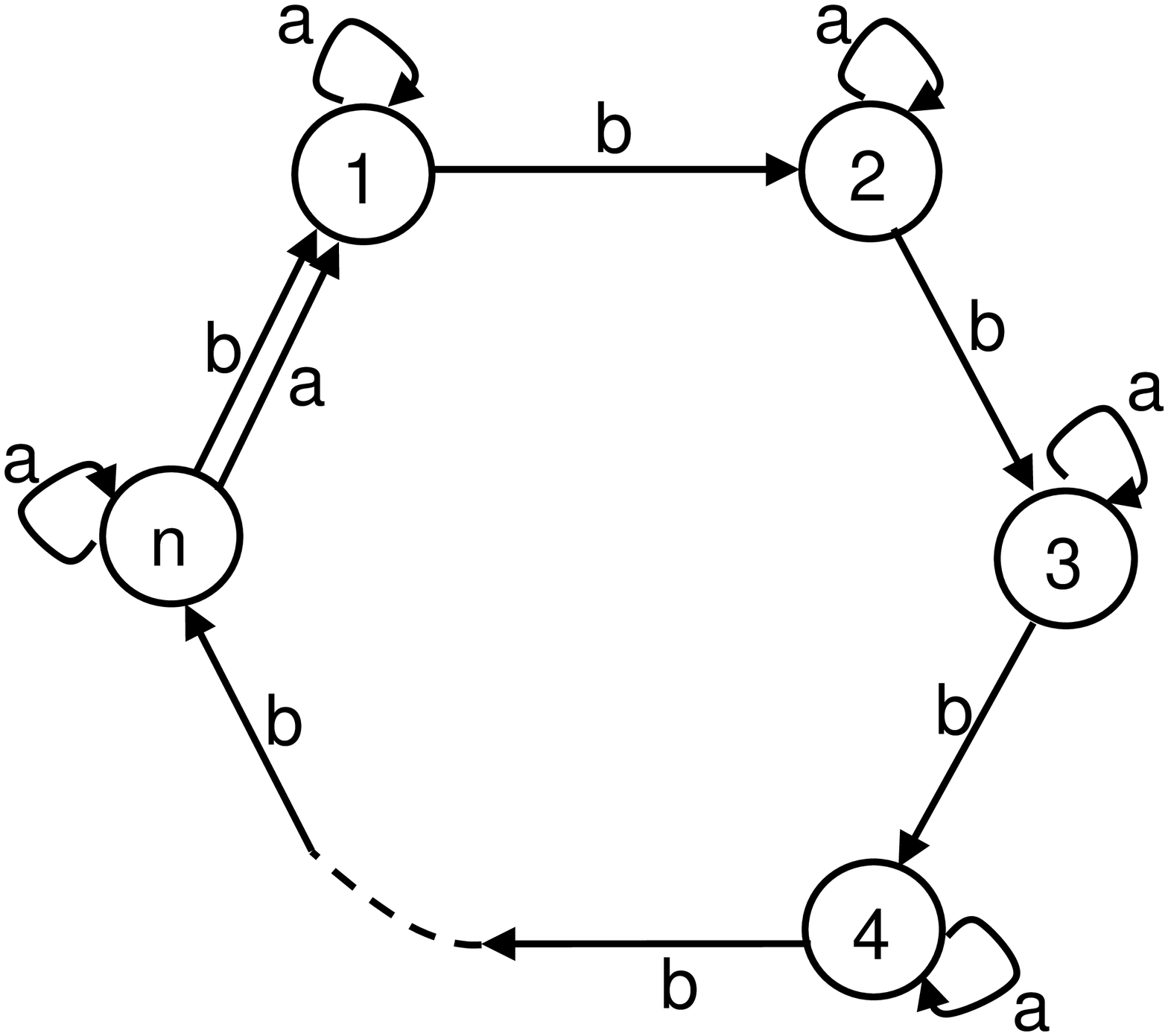} \\(b) 
\end{tabular} 
\caption{Construction of our extremal example. (a) \cerny's well known automaton, whose shortest synchronizing words are of length $(n-1)^2$.  This set of matrices is not primitive. (b) Our extremal example, which is primitive.  It has no positive product of length shorter than $( n^2/2).$  }\label{fig-extremal}\end{figure}

\begin{conj}\textbf{\v Cern{\'y}'s conjecture, 1964} \label{conj-cerny-initial}\cite{cernyPirickaRosenauerova64} 
Let $\cm \subset \{0,1\}^{n\times n}$ be a synchronizing automaton.  Then, there is a synchronizing word of length at most $(n-1)^2.$\end{conj} 
In fact, it is even not known whether there is a valid bound with a quadratic growth in $n,$ and we study in the rest of this paper the weaker following conjecture. 
\begin{conj}\label{conj-cerny} 
Let $\cm \subset \{0,1\}^{n\times n}$ be a synchronizing automaton.  Then, there is a synchronizing word of length at most $Kn^2$ for some fixed $K>0.$\end{conj} 
Fig. \ref{fig-extremal} (a)  represents a synchronizing automaton whose shorter synchronizing word is of length $(n-1)^2,$ as proved in \cite{cernyPirickaRosenauerova64}.  Thus, if Conjecture \ref{conj-cerny-initial} is true, the bound in the conjecture is tight. 
  
We first present a technical result which makes the bridge between the combinatorial problem studied in the present paper and the notion of synchronizing automaton.   
 
\begin{thm} \label{prop-prim-sync} 
For any primitive set of nonnegative matrices $$\setmat=\{A_1,\dots,A_m\} \subset \{0,1\}^{n\times n}$$ satisfying Assumption \ref{assum-nzrnzc} there exists a synchronizing automaton $$\setmat'=\{A'_1,\dots,A'_t\}$$ such that $$\forall 1\leq s\leq t,\, \exists l\in \{1,\dots,m\}: A'_s\leq A_l \, {\mbox{(entrywise)}}. $$ 
\end{thm} 
We attract the attention of the reader to the fact that the number of matrices in the automaton is not necessarily the same as the number of matrices in the initial set $\setmat.$ 
\vskip -1.5pc \begin{pf}  Let us consider the positive product $A_{i_1}A_{i_2}\dots A_{i_t}\in \cm^*.$ We will keep the different matrices $A_{i_l}:$ $l=1\dots t$ for the construction of our automaton.  Since the product is positive, there are actually paths from nodes $1,\dots,n$ to node (say,) $1,$ in the sequence of graphs $G_{i_1}G_{i_2}\dots G_{i_t}.$ In order to obtain our automaton, we have to remove edges (i.e., put some entries to zero in our constructed matrices $A_i'$) so that one and only one entry in every row is equal to one.  If we manage to do that while keeping the $n$ paths in the sequence of our corresponding graphs $G'_{i_1}G'_{i_2}\dots G'_{i_t},$ then we will have a synchronizing automaton.   

In order to do that, we simply keep in each matrix $A_{i_l}$ all the edges that are part of the above mentioned paths.  If there is a node $v$ in the graph $G_{i_l}$ such that $v$ is on none of these paths (at level $l$), we can just pick any edge leaving $v$ in order to define a valid automaton.  Such an edge exists because all matrices in $\setmat$ have nonzero rows and columns.\\ 
Now, there might be some graphs $G_{i_l}$ in which two edges leaving the same node $v$ have been kept.  However, this could occur only if there are two separate paths leaving $v$ (at the level $l$) and reaching node $1$ at level $t$.  Thus, one can safely iteratively remove these edges in excess, and making sure at the same time that if a node was connected to the node 1 at level $t,$ it remains connected by at least one path throughout this process. \qedsymbol \end{pf} 
 
\begin{thm}\label{thm:bound-sync-2bound-prim} For any primitive set of matrices $\cm$ of dimension $n$ satisfying Assumption \ref{assum-nzrnzc} there is a product of length smaller than $2f(n)+n-1$ with positive entries, where $f(n)$ is any upper bound on the minimal length of a synchronizing word for $n$-dimensional automata. 
\end{thm} 
\vskip -1.5pc \begin{pf} 
From Theorem \ref{prop-prim-sync} above, let us consider the automaton $A'$ whose matrices are smaller (entrywise) than matrices from $\cm.$  There is a product $B_1$ of length $f(n)$ with a positive column (say, the $i$th one).  Now, reasoning on the set $\cm^{Tr},$ there is a product $B_2$ of length $f(n)$ with a positive row (say, the $j$th one).  Now one can take a product $C$ (of length smaller than $n$) such that $C_{i,j}>0,$ and one obtains $$B_1CB_2>0.$$ \qedsymbol
\end{pf} 
  
\begin{cor} 
For any primitive set of matrices of dimension $n$ satisfying Assumption \ref{assum-nzrnzc} there is a product of length smaller than $$(n(7n^2+6n+8)-24)/24=O(n^3) $$ with positive entries. 
\end{cor} 
\vskip -0.3pc \begin{pf} It is known \cite{trahtman2011} that any synchronizing automaton has a synchronizing word of length smaller or equal to 
$$f(n)=n(7n^2+6n-16)/48 .$$  By combining this bound with Theorem \ref{thm:bound-sync-2bound-prim} above, we obtain the result. \qedsymbol
\end{pf} 

We did not try to optimize the bound in the above theorem.  Most probably simple arguments could allow to lower it with the same general ideas. \rmj{It also shows that the upper bound given in \cite{voynov} is not sharp. It is natural to ask whether a cubic upper bound is sharp; this is problem 
1 in \cite{voynov}.  One might further }  hope to decrease the bound to a quadratic degree, as formalized in the next conjecture:

\begin{conj} \label{conj-cerny-equiv} 
There is a constant $K$ such that for any set of primitive matrices of dimension $n$ satisfying Assumption \ref{assum-nzrnzc} there is a product of length smaller than $Kn^2$ with positive entries. 
\end{conj} 
\rmj{Note that Conjecture \ref{conj-cerny} being true would directly imply Conjecture \ref{conj-cerny-equiv}.}
We finish by providing a lower bound for the shortest length of a positive product. 

\begin{ex}\label{ex-extremal} 
Fig. \ref{fig-extremal} (b) represents 
a set of matrices which is primitive, but the length of any positive product is at least $n^2/2.$ To see this, let us consider only the first row of a product of length $t,$ which we denote $v_t,$ and start with the empty product (i.e. the Identity matrix).  Let us denote $A_a,A_b$ the two matrices corresponding to the graph in Fig. \ref{fig-extremal} (b).  Observe that right-multiplying our product with $A_b$ shifts all the entries to the right, that is, 
$$  v_tA_b (i) =v_t(j):j=i-1 \mod n. $$ Also, multiplying our product with $A_a$ leaves $v_t$ unchanged, except if $v_t(1)=0,v_t(n)>0, $ in which case $v_{t+1}(1)=v_t(n).$ Thus, it is straightforward to prove inductively that the only way to increase the number of nonzero entries in a vector $v_t=(1,0,\dots,0,1,\dots,1)$ is to apply $A_b$ $n-1$ times (or $k(n-1)$ for some natural number $k$), followed by $A_a$ in which case the vector has the same general shape $v_t=(1,0,\dots,0,1,\dots,1), $ with one more $1$.  Thus, any positive product having all its entries in the first row positive, this process has to be repeated $n-1$ times, which brings a lower bound of $n(n-1).$  
\end{ex} 
%

As a direct consequence of Example \ref{ex-extremal}, we have the following corollary.
\begin{cor} 
The upper bound in Conjecture \ref{conj-cerny-equiv} cannot be $o(n^2).$ 
\end{cor}

%

\bibliographystyle{plain} 
\bibliography{references-primitive} 
%
%
%
%
%
%
%
\end{document}